\documentclass[12pt, reqno]{amsart}
\usepackage{amssymb}
\usepackage{amsmath}
\usepackage[utf8]{inputenc}
\usepackage{longtable}
\usepackage{cite}
\usepackage[all,cmtip]{xy}
\usepackage{color}
\usepackage{fullpage}

\usepackage{amsmath}
\usepackage{amsfonts}
\usepackage{amssymb}
\usepackage{amsthm}
\usepackage{mathtools}
\usepackage{mathrsfs}
\usepackage{float}
\usepackage[all,cmtip]{xy}
      \theoremstyle{definition}
\newtheorem{defi}{Definition}[section]
\theoremstyle{plain}
\newtheorem{thm}[defi]{Theorem}

\newtheorem{prop}[defi]{Proposition}
\newtheorem{cor}[defi]{Corollary}
\newtheorem{lemma}[defi]{Lemma}
\theoremstyle{remark}

\newtheorem{rmk}[defi]{Remark}

\theoremstyle{definition}

\newtheorem*{ack}{Acknowledgement}

\newcommand{\ra}{\rightarrow}

\newcommand{\longra}{\longrightarrow}

\newcommand{\ka}{{\mathcal A}}

\newcommand{\km}{{\mathcal M}}

\newcommand{\ko}{{\mathcal O}}

\newcommand{\IN}{{\mathbb N}}

\newcommand{\Q}{{\mathbb Q}}
\newcommand{\R}{{\mathbb R}}

\newcommand{\Z}{{\mathbb Z}}

\newcommand{\gp}{{\mathfrak p}}

\newcommand{\rH}{{\rm H}}

\newcommand{\im}{\operatorname{im}}

\newcommand{\Spec}{\operatorname{Spec}}

\newcommand{\Pic}{\operatorname{Pic}}

\newcommand{\End}{\operatorname{End}}

\newcommand{\Gal}{\operatorname{Gal}}

\newcommand{\Hom}{\operatorname{Hom}}

\newcommand{\NS}{\operatorname{NS}}

  \makeatletter
\newcommand{\xdashrightarrow}[2][]{\ext@arrow 0359\rightarrowfill@@{#1}{#2}}
\newcommand{\xdashleftarrow}[2][]{\ext@arrow 3095\leftarrowfill@@{#1}{#2}}
\newcommand{\xdashleftrightarrow}[2][]{\ext@arrow 3359\leftrightarrowfill@@{#1}{#2}}
\def\rightarrowfill@@{\arrowfill@@\relax\relbar\rightarrow}
\def\leftarrowfill@@{\arrowfill@@\leftarrow\relbar\relax}
\def\leftrightarrowfill@@{\arrowfill@@\leftarrow\relbar\rightarrow}
\def\arrowfill@@#1#2#3#4{%
  $\m@th\thickmuskip0mu\medmuskip\thickmuskip\thinmuskip\thickmuskip
   \relax#4#1
   \xleaders\hbox{$#4#2$}\hfill
   #3$%
}
\makeatother

	\title{On the Picard numbers of abelian varieties}
	\author{Klaus Hulek \and Roberto Laface}
	
	\address{Klaus Hulek \newline Institut f\"{u}r Algebraische Geometrie, Leibniz Universit\"{a}t Hannover, Welfengarten 1 30167 Hannover (Germany)
}
    \email{hulek@math.uni-hannover.de}
	\address{Roberto Laface \newline Technische Universit\"at M\"unchen, Zentrum Mathematik - M11,
    Boltzmannstra{\ss}e 3,
    D-85478 Garching bei M\"unchen (Germany)
}
    \email{laface@ma.tum.de}

\begin{document}

\thispagestyle{empty}
\begin{abstract}
In this paper we study the possible Picard numbers $\rho$ of an abelian variety $A$ of dimension $g$. It is well known that this satisfies the inequality $1 \leq \rho \leq g^2$. 
We prove that the set $R_g$ of realizable Picard numbers of abelian varieties of dimension $g$ is not complete for every $g \geq 3$, namely that $R_g \subsetneq [1,g^2] \cap \IN$. Moreover, we study the structure of $R_g$ as $g \ra +\infty$, and from that we deduce a structure theorem for abelian varieties of large Picard number. In contrast  to the non-completeness of any of the sets $R_g$ for $g \geq 3$, we also show that the Picard numbers of abelian varieties are asymptotically complete, i.e. $\lim_{g \ra +\infty} \#R_g / g^2 = 1$. As a byproduct, we deduce a structure theorem for abelian varieties of large Picard number.
Finally we show that all realizable  Picard numbers in $R_g$ can be obtained by an abelian variety defined over a number field.
\end{abstract}
	\maketitle
\setcounter{tocdepth}{1}

\section{Introduction}

For an algebraic variety $X$ over the field of complex numbers the Lefschetz $(1,1)$-theorem says that the N\'eron-Severi group 
\[ \NS(X) = \rH^2(X,\Z) \cap \rH^{1,1}(X).\]
Consequently, the rank $\rho(X)$ of the N\'eron-Severi group, the so-called {\em Picard number}, satisfies the inequality $1 \leq \rho(X) \leq h^{1,1}(X)$. 
Computing the Picard number is in general a difficult question, as already the case of projective surfaces shows. For example, the Picard number of a quintic surface
$S$ in $\mathbb P^3$ satisfies the inequality $\rho(S) \leq 45$. It is known that all  numbers between $1$ and $45$  can be obtained if 
one allows the surface to have $ADE$-singularities, but it remains an open problem for smooth surfaces, where the maximum known is $41$ \cite{schuett11}, \cite{schuett15}.

In this article we will concentrate on the Picard numbers of abelian varieties.  
To put this into perspective it is worthwhile to recall  the situation for surfaces. For abelian surfaces all possible Picard numbers between $1$ (or $0$ if one includes the
non-algebraic case) and $4$ occur.  Indeed, a very general abelian surface has $\rho=1$, whereas Picard numbers from 2 to 4 can be realized by taking a product $E_1 \times E_2$ of two elliptic curves. If the two elliptic curves are not isogenous, then $\rho=2$, if they are isogenous but they do not have complex multiplication, then $\rho =3$, while if they also have complex multiplication $\rho =4$. For the other surfaces with trivial canonical bundle the situation is similar: for K3 surfaces 
all possibilities between $1$ (respectively $0$) and $20$ can occur, as can be seen by the Torelli theorem for K3 surfaces. Enriques surfaces and bi-elliptic surfaces have no holomorphic 2-forms, and thus their Picard number is $10$ and $2$ respectively.

For higher-dimensional varieties with numerically trivial canonical bundle the situation is as follows.
By the Beauville-Bogomolov decomposition theorem \cite{beauville84}, every K\"ahler manifold with trivial first Chern class admits a finite cover which is a product of tori,
Calabi-Yau varieties and irreducible holomorphic symplectic manifolds (IHSM), also know as hyperk\"ahler ma\-ni\-folds. For higher dimensional Calabi-Yau varieties $Y$ we always 
have $\rho(Y)=b_2(Y)$ as $h^{2,0}(Y)=h^{0,2}(Y)=0$. For  irreducible holomorphic symplectic manifolds $X$ one can use Huybrechts' surjectivity of the period map \cite{huybrechts99} to conclude, as in the case of K3 surfaces, that all values $0 \leq \rho(X) \leq b_2(X) - 2$ can be obtained. This leaves us with the case of abelian varieties which is the topic of this note. Surprisingly little seems to be known about the possible Picard numbers of abelian varieties. Our aim is to make a first start to remedy this situation, using mostly elementary methods.

Let $A$ be a complex torus of dimension $g$. Its cohomology is the exterior algebra over $\rH^1(A,\mathbb Z) \cong {\mathbb Z}^{2g}$. 
In particular, this implies that the $k$\textsuperscript{th} Betti numbers are $b_k(A) = {2g \choose k}$.  As $\rH^{p,0}(A) \cong \rH^0(A, \Omega_A^p)$, we get $h^{p,0}(A) = {g \choose p}$, and thus $h^{1,1}(A) =  g^2$. 
We shall from now on exclude the case of non-algebraic tori and concentrate on abelian varieties.
By the above we  know that 
\[1 \leq \rho \leq g^2.\]
As we have already seen, any number $1 \leq \rho(A) \leq 4$ can be achieved for abelian surfaces. However, the situation changes significantly in higher dimension. 

Given an arbitrary abelian variety $A$, we can invoke Poincar\'e Complete Reducibility Theorem \cite[Thm. 5.3.7]{birkenhake-lange04} to pass to a better representative in its isogeny class, namely
\[A \longra A_1^{n_1} \times \cdots \times A_r^{n_r},\]
where $A_i$ is a simple abelian variety ($i =1, \dots, r$), and $A_i$ is not isogenous to $A_j$ if $i \neq j$. Moreover, the abelian varieties $A_i$ and the integers $n_i$ are uniquely determined up to isogeny and permutations. A result of Murty \cite[Lemma 3.3]{murty84} describes the Picard number of a self-product $B^k$ of a simple abelian variety in terms of $k$ and the dimension of $B$. In light of this, we prove in Proposition \ref{pic} a splitting result concerning the Picard group of varieties of the form $A \times B$ with $\Hom(A,B)=0$, which allows us to compute the Picard number of such products. This, together with Murty's result in \cite{murty84}, provides us with a theoretical algorithm for computing the Picard number of a given abelian variety. 

It is then a combinatorial question as to determine the set $R_g$ of possible Picard numbers of abelian varieties for a given genus $g$. Very little seems to be known about this. The purpose of this paper is to take a first step in the analysis of $R_g$. 
As a first result we show that there are gap series for the possible Picard numbers of abelian varieties, and therefore that the sets $R_g$ are not \textit{complete} for every $g \geq 3$. In fact, it is not hard to show that $R_3= \lbrace 1, \dots, 6,9 \rbrace$: indeed, a very general abelian threefold has Picard number $\rho=1$; a product $S \times E$, $S$ being a very general abelian surface and $E$ being an elliptic curve, has Picard number $\rho=2$; all other Picard numbers can be obtained by using suitable products of elliptic curves. This phenomenon had previously been notice by Shioda in \cite[Appendix]{shioda75}.

As the dimension $g$ grows larger, clear gaps in the set of possible Picard numbers start to appear. Moreover, more and more gaps occur as $g \ra \infty$. The following result shows the existence of two precise gaps and characterizes the three largest Picard numbers for an abelian variety.

\begin{thm}\label{mainthm1}
\noindent
\begin{itemize}
\item[(1)] Fix $g \geq 4$. There does not exist any abelian variety of dimension $g$ with Picard number $\rho$ in the following range:
\[ (g-1)^2+1 < \rho < g^2.\]
\item[(2)] Fix $g \geq 7$. There does not exist any abelian variety of dimension $g$ with Picard number $\rho$ in the following range: 
\[(g-2)^2 + 4 < \rho < (g-1)^2 +1.\]
\end{itemize}
\end{thm}

We would like to remark that the conditions on the dimension $g$ given in Part 1 and 2 of Theorem \ref{mainthm1} are necessary. In fact, as for Part 1, for $g=2$ all Picard numbers occur, and for $g=3$ there exists an abelian threefold of Picard number $\rho=6$ (namely, the product of three isogenous elliptic curves without complex multiplication). Similar considerations can be made for Part 2 of Theorem \ref{mainthm1} and $g\leq 6$. 
After some preliminary work in Section \ref{sec:prelim}, we shall prove this theorem in Section \ref{sec:proofmain}.
As an application of our analysis  we derive in Section \ref{sec:Structure} a structure theorem for abelian varieties with large Picard number, namely Theorem \ref{structurethm}. 

The above results are a first indication of a much more general phenomenon which we study more systematically in Section \ref{sec:asymptotics}, where we consider the 
behaviour of the set $R_g$ asymptotically, namely as $g$ grows.
In particular, we define the \textit{asymptotic density} of Picard numbers to be the quantity 
\[ \delta := \lim_{g \ra + \infty} \frac{\# R_g}{g^2}. \]
Contrary to the non-completeness of any of the sets $R_g$, we prove asymptotic completeness in Theorem \ref{density}, namely 

\begin{thm}\label{mainthm2} 
The Picard numbers of abelian varieties are \textit{asymptotically complete}:
\[ \delta = \lim_{g \ra + \infty} \frac{\# R_g}{g^2}=1. \]
\end{thm}

By using similar techniques, we are also able to describe the distribution of the Picard numbers within $[1,g^2] \cap \IN$ (Theorem \ref{asymptoticstr}). As a consequence, we obtain a structure theorem for abelian varieties of large Picard number in Corollary \ref{cor:morestructure}. 
We also provide a practical algorithm which allows to compute the sets $R_g$ inductively. 

Finally we show that all realizable Picard numbers $\rho \in R_g$ can be obtained by an abelian variety defined over a number field.
 
\begin{ack}
We would like to thank Bert van Geemen for his genuine interest in this question, and Matthias Sch\"utt for discussions around this topic and beyond. RL would like to particularly thank Fran\c{c}ois Charles for the many fruitful discussions and for his invitation to the IH\'ES, whose excellent working conditions are here gratefully acknowledged. We would like to thank Davide Lombardo for suggesting the proof of Lemma \ref{lombardo}. 
We gratefully acknowledge a very fruitful exchange of ideas with Ben Moonen who also supplied the proof of Theorem \ref{teo:endonumberfield}. Finally, we are grateful to the anonymous referee for his or her comments, which have helped improving the manuscript. 
This research was partially funded by the DFG funded GRK 1463 "Analysis, Geometry and String Theory".
\end{ack}

\section{Preliminary work}\label{sec:prelim}

In this section we will develop the basic tools of our analysis. Some of these results are of independent interest in their own right.

\subsection{Additivity of the Picard number for non-isogeneous products}
As the Picard number of an abelian variety is invariant under isogenies \cite[Ch.\,1, Prop.\,3.2]{birkenhake-lange99}, we can pick a convenient representative in its isogeny class. Such a choice is indicated by the following result \cite[Thm 5.3.7]{birkenhake-lange04}:

\begin{thm}[Poincar\'e's Complete Reducibility Theorem]\label{poincare}
	Given an abelian variety $A$, there exists an isogeny
	\[A \longra A_1^{n_1} \times \cdots \times A_r^{n_r},\]
	where $A_i$ is a simple abelian variety ($i =1, \dots, r$), and $A_i$ is not isogenous to $A_j$ if $i \neq j$. Moreover, the abelian varieties $A_i$ and the integers $n_i$ are uniquely determined up to isogeny and permutations.
\end{thm}

Let us now consider a product of simple abelian varieties as in Theorem \ref{poincare}. The fact that $A_i $ is not isogenous to $A_j$ for $i \neq j$ yields the following  splitting of the Picard group:

\begin{prop}\label{pic}
	Let $A_1, \dots, A_r$ be simple abelian varieties, such that $A_i$ is not isogenous to $A_j$ for $i \neq j$. Then the (exterior) pullback of line bundles yields an isomorphism
	\[ \prod_{i=1}^r \Pic(A_i^{n_i}) \cong \Pic \Bigg(\prod_{i=1}^r A_i^{n_i} \Bigg).\]
\end{prop}

Clearly, exterior pull-back of line bundles always yields an injective map, but surjectivity is a special feature. In fact, if $E$ is an elliptic curve, the abelian surface $E \times E$ has Picard number $\rho \in \lbrace 3,4 \rbrace$, depending on the presence of CM. Therefore, the exterior pull-back map
\[ \Pic(E) \times \Pic(E) \longra \Pic(E \times E)\]
cannot be surjective, as otherwise we would get a surjective map of the corresponding N\'eron-Severi groups, hence yielding a contradiction, since $\NS(E) \cong \Z$.

\begin{proof}
	Exterior pullback of line bundles 
	$$
	\psi(L_1, \ldots, L_r)=  L_1  \boxtimes  \cdots  \boxtimes  L_r
	$$
	defines the following commutative diagram
	\begin{equation*}
	\xymatrix{
		0   \ar[r]  & \Pic^0(\prod_{i=1}^r  A_i^{n_i})   \ar[r]  & \Pic(\prod_{i=1}^r  A_i^{n_i}) \ar[r] & \NS(\prod_{i=1}^r  A_i^{n_i}) \ar[r] & 0  & \\
		0 \ar[r] \ar[u] & \prod_{i=1}^r \Pic^0(A_i^{n_i})  \ar[u]_{\psi^0} \ar[r]& \prod_{i=1}^r \Pic(A_i^{n_i})  \ar[u]_{\psi} \ar[r]& \prod_{i=1}^r \NS(A_i^{n_i}) \ar[u]_{\psi^{\NS}} \ar[r] & 0 .\ar[u]\\
	}
	\end{equation*}
	We will show that $\psi^0$ and $\psi^{\NS}$ are isomorphisms, thus proving the proposition. Clearly $\psi^0$ is injective, and since $\psi^0$ is a homomorphism of abelian varieties of the 
	same dimension it must be an isomorphism.  To prove that $\psi^{\NS}$ is an isomorphism we recall from \cite[Ch.~2]{birkenhake-lange04} that a polarization on an abelian variety 
	$A$ is  given by a finite isogeny $f: A \to A^\vee$ whose analytic representation is hermitian.  By assumption the abelian varieties $A_i$ and $A_j$ are not isogeneous for $i \neq j$.
	Hence $\Hom(A_i,A_j)=\Hom(A_i,A_j^\vee)=0$ and every isogeny 
	\[f: \, \prod_{i=1}^r  A_i^{n_i} \longra (\prod_{i=1}^r  A_i^{n_i})^\vee\]
	is of the form $f=(f_1, \ldots, f_r)$ where 
	$f_i: A_i^{n_i} \longra (A_i^{n_i})^\vee$ is an isogeny. Since a direct sum of endomorphisms is hermitian if and only if all its summands are, the claim follows, as every class in the N\'eron-Severi group can be written as the difference of two ample classes (i.e.~two polarizations). 
\end{proof}

As a consequence, we get that the Picard number is additive (but not strongly additive) for product varieties coming from the Poincar\'e's Complete Reducibility Theorem.

\begin{cor}\label{picnum}
	Let $A_1, \dots, A_r$ be simple abelian varieties, such that $A_i$ is not isogenous to $A_j$ for $i \neq j$. Then,
	\[ \rho \Bigg( \prod_{i=1}^r A_i^{n_i} \Bigg) = \sum_{i=1}^r \rho (A_i^{n_i}).\]
\end{cor}

\subsection{Picard numbers of self-products}\label{subsec:selfproducts}

Due to additivity, we are left to see how to compute the Picard number in the case of a self-product of a simple abelian variety. 
The endomorphism ring $\End(A)$ of an abelian variety $A$ is a finitely generated free abelian group and 
hence $F:= \End(A)\otimes \Q \equiv \End_\Q(A)$ is a finite dimensional $\Q$-algebra.
Any polarization $L$ on $A$ defines an involution on $F$, the so-called Rosati involution.  
If $A$ is a simple abelian variety, then $F$ is a finite-dimensional skew field admitting  a positive anti-involution. 
Such pairs were classified by
Albert \cite{albert34}, \cite{albert35}, see also \cite[Proposition 5.5.7]{birkenhake-lange04} for a summary.

Let $K$ be the centre of $F$. We will say that $F$ is \textit{of the first kind} if the Rosati involution acts trivially on $K$, and \textit{of the second kind} otherwise. Let us denote by $K_0$ the maximal real subfield of $K$, and let us consider the following invariants of $F$:
\[ [F:K] =d^2, \qquad [K:\Q] = e, \qquad [K_0:\Q] = e_0. \]
Notice that, as $K$ is the center of $F$, $[F:K]$ is always a square. 

As a useful example, we can consider quaternion algebras over a field $K$. Let us recall the reader that any quaternion algebra $F/K$ is either a division ring or it is isomorphic to $M_2(K)$. In light of this, we can define the \textit{ramification locus} of $F$ as 
\[\text{Ram}(F) = \lbrace \gp \in \Spec \ko_K \, \vert \, \text{$K_\gp$ is a division ring} \rbrace.\]
A quaternion algebra is ramified at a finite even number of places. The ramification locus governs the isomorphism classes of quaternion algebras: there is a 1:1 correspondence between isomorphism classes of quaternion algebras and subsets of non-complex places of $K$ of even cardinality. 

By \cite[Proposition 5.5.7]{birkenhake-lange04} the classification divides into four types, where the first three are of the first kind:

\begin{enumerate}
	\item \textit{Type I}: $F$ is a totally real number field, so that $d=1$ and $e =e_0$.
	\item \textit{Type II}: $F$ is a totally indefinite quaternion algebra over a totally real number field $K$, 
	i.e. 
	\[ \emptyset = \text{Ram}(F) \cap \lbrace \text{archimedean places of $K$} \rbrace \quad \text{and} \quad \text{Ram}(F) \neq \emptyset. \]
	In particular, we have that $d=2$ and $e = e_0$.
	\item \textit{Type III}: $F$ is a totally definite quaternion algebra over $K$, i.e.
	\[ \emptyset \neq \text{Ram}(F) \supseteq \lbrace \text{archimedean places of $K$} \rbrace \]
	holds. Again, we have $d=2$ and $e =e_0$.
	\item \textit{Type IV}: $F$ is of the second kind, and it center $Z(F)=K$ is a CM field with maximal real subfield $K_0$.
\end{enumerate}

The following result, due to Murty, gives a complete description of the Picard number of a self-product of a simple abelian variety.

\begin{prop}[Lemma 3.3 of \cite{murty84}]\label{murty}
	Let $A$ be a simple abelian variety. Set $e:=[K : \Q]$, $d^2 :=[F:K]$. Then, for $k \geq 1$, one has
	\[
	\rho(A^k) = 
	\begin{cases} 
	\frac{1}{2}ek(k+1) & \text{Type I}\\
	ek(2k+1) & \text{Type II}\\
	ek(2k-1) & \text{Type III}\\
	\frac{1}{2}ed^2k^2 & \text{Type IV.}
	\end{cases}
	\]
\end{prop}

In fact, Murty's result is in terms of the maximal commutative subfield $E$ of $F$, which has degree $[E:K]$ over $K$. However, the proof of \cite[Lemma 2.2]{murty84} implies that $[E:K]^2 = [F:K]$. Proposition \ref{murty} enables us to compute the following bound for the Picard number of a self-product of a simple abelian variety:

\begin{cor}\label{lemma1}
	Let $A$ be a simple abelian variety of dimension $n$, and let $k \geq 1$. Then $\rho(A^k) \leq \frac{1}{2}nk(2k+1)$.
\end{cor}

\begin{proof}
	Proposition \ref{murty} applied with $k=1$ allows us to compute the Picard number of $A$:
	\[\rho = \rho(A) = 
	\begin{cases} 
	e & \textit{Type I}\\
	3e & \textit{Type II}\\
	e & \textit{Type III}\\
	\frac{1}{2}ed^2 & \textit{Type IV.}
	\end{cases}
	\]  
	Now, plugging this back in Proposition \ref{murty} gives the following  reformulation in terms of the Picard number of $A$:
	\[ \rho(A^k) = 
	\begin{cases}
	\frac{1}{2}\rho k(k+1) & \textit{Type I}\\
	\frac{1}{3}\rho k(2k+1) & \textit{Type II}\\
	\rho k (2k-1) & \textit{Type III}\\
	\rho k^2 & \textit{Type IV.}
	\end{cases}
	\]
	The divisibility conditions for  $\rho$ given by  \cite[Prop.~5.5.7]{birkenhake-lange04} imply that
	\[ \rho \leq 
	\begin{cases} 
	n & \textit{Type I}\\
	\frac{3}{2}n & \textit{Type II}\\
	\frac{1}{2}n & \textit{Type III}\\
	n & \textit{Type IV}
	\end{cases}
	\]
	and, therefore, we see that
	\[ \rho(A^k) \leq 
	\begin{cases}
	\frac{1}{2}n k(k+1) & \textit{Type I}\\
	\frac{1}{2}n k(2k+1) & \textit{Type II}\\
	\frac{1}{2}n k (2k-1) & \textit{Type III}\\
	n k^2 & \textit{Type IV}
	\end{cases}
	\]
	from which the result follows.
\end{proof}

In the case of a self-product of an elliptic curve this gives the well known 

\begin{cor}\label{cor:ellipticproducts} 
If $E$ is an elliptic curve, then
\[ \rho(E^k) = 
	\begin{cases}
	\frac{1}{2} k(k+1) & \textit{$E$ has no CM}\\
	k^2 & \textit{$E$ has CM.}
	\end{cases}
	\] 

\end{cor}

We will use these results, in particularly the case of a self-product of elliptic curves,  frequently in the proof of Theorem \ref{mainthm1}. Notice that Corollary \ref{lemma1} provides us with a bound on the Picard number of $A^k$ which is independent of the type of the endomorphism ring of $A$.

\section{Restrictions on the Picard number}\label{sec:proofmain}

\subsection{Some bounds on the Picard number}
We would like to show that there are better bounds on the Picard number, if one is given a partition of the dimension. More precisely, letting $A$ be an abelian variety, we define $r(A)$ to be the \textit{length} of a decomposition according to Poincar\'e Complete Reducibility Theorem. In other words, given an abelian variety $A$, Theorem \ref{poincare} gives an isogeny 
\[A \longra A_1^{n_1} \times \cdots \times A_r^{n_r}, \]
and we set $r(A) := r$. Notice that this quantity is well-defined because the factors $A_i$ and the powers $n_i$ are determined up to permutations and isogenies. Then, for $r \leq g$, we define $M_{r,g}$ as
\[ M_{r,g} := \max \lbrace \rho(A) \, \vert \, \dim A =g, \ r(A) = r \rbrace .\]
In other words, $M_{r,g}$ is the largest Picard number that can be realized by a $g$-dimensional abelian variety that splits into a product of $r$ non-isogenous pieces in its isogeny class.

\begin{prop}\label{bounds}
	For integers $r, g \in \IN$ such that $r \leq g$, one has $M_{r,g} = [ g -(r-1)]^2 + (r-1)$. This value is attained as the Picard number of $E^{g-r+1} \times E_1 \times \cdots \times E_{r-1}$, where $E$ is a CM elliptic curve not isogenous to any of the $E_i$'s, and $E_i$ and $E_j$ are not isogenous for $i \neq j$.
\end{prop}

\begin{proof}
	If $A \sim A_1 \times \cdots \times A_r$, $\Hom(A_i, A_j)=0$ for $i \neq j$, then
	\[ \rho(A) \leq k_1^2 + \cdots + k_r^2\]
	where $k_i := \dim A_i$ ($i=1, \dots, r$) and $k_1 + \cdots + k_r = g$. Hence we are looking for the maxima of the function 
	\[ h(x_1, \dots, x_r) = x_1^2 + \cdots + x_{r-1}^2 + x_r^2 \] 
	on the integral points of the  simplex 
	\[
	\Omega_{r,g}=\{(x_1, \ldots, x_r) \mid x_i \geq 1,  x_1+ \ldots x_r=g\}. 
	\]
	These points are precisely the vertices 
	\[ \big\lbrace (g-r+1, 1 ,\dots, 1), (1, g-r+1, 1,\dots, 1), \dots,  (1, \dots, 1, g-r+1) \big\rbrace. \]
	By the symmetry of $h$, the maximum is attained at any of these points, with value 
	\[h(g-r+1, 1, \dots,1)=[ g- (r-1)]^2 + (r-1).\]
	Therefore, we conclude that $\rho(A) \leq [ g- (r-1)]^2 + (r-1)$. By applying Proposition \ref{pic} and Corollary \ref{cor:ellipticproducts} one can see that the abelian variety 
	\[ E^{g -r+1} \times E_1 \times \cdots \times E_{r-1}\]
	with $E$ a CM curve and $E,E_i,E_j$ for $i\neq j$ not pairwise mutually isogeneous,
	has Picard number $[ g -(r-1)]^2 + (r-1)$, and we are thus done. 
\end{proof}

\begin{cor}\label{lemmastr}
	Let $A$ be an abelian variety. Then, 
	\[ \rho(A)= M_{r(A),g} \Longleftrightarrow A \sim E^{g-(r-1)} \times E_1 \times \cdots \times E_{r-1}, \]
	where $E$ is a CM elliptic curve not isogenous to any of the $E_i$'s, and $E_i$ and $E_j$ are not isogenous for $i \neq j$.
\end{cor}

\begin{rmk}
	The numbers $M_{r,g}$ give the following (strictly) increasing sequence of positive integers:
	\[ g=M_{g,g} < M_{g-1,g} < \cdots < M_{3,g} < M_{2,g} < M_{1,g}=g^2. \]
\end{rmk}

We will now proceed with the proof of Theorem \ref{mainthm1}, which we divide into two parts.
\subsection{Proof of part (1)}
Let $A$ be an abelian variety of dimension $g \geq 4$ with Picard number $\rho = \rho(A)$. We will divide our analysis of the Picard number $\rho$ into the following cases:
\begin{enumerate}
	\item[(a)] $A$ has length at least two, i.e.~$r(A)\geq 2$;
	\item[(b)] $A$ is a self-product of a lower dimensional abelian variety.
\end{enumerate}

\textbf{Case (a).} Since $r(A) \geq 2$, we have that $A \sim A_1 \times A_2$ with $\Hom(A_1,A_2) = 0$. Let $n :=\dim A_1$, so that $\dim A_2 = g-n$. Then, $\rho(A) \leq n^2 + (g-n)^2$. Consider the function 
\[f(x):=x^2 + (g-x)^2\]
on $\Omega =[1,g-1]$. It attains its maximum at $x=1$ and $x=g-1$, with value $f(1)=f(g-1)=(g-1)^2 +1$. Therefore, $\rho(A) \leq (g-1)^2+1$.

\textbf{Case (b).} 	Let $B$ be an $m$-dimensional simple abelian variety, and suppose $A$ is isogenous to $B^k$, for $k := g/m$. If $m =1$ (i.e.\,$B$ is an elliptic curve), then
again by Corollary \ref{cor:ellipticproducts}
\[ \rho(B^g) =
\begin{cases}
{{g+1} \choose 2} & \text{$B$ has no CM}\\
g^2 & \text{$B$ has CM.}
\end{cases}
\]

If $B$ has CM, then $A$ attains the maximal Picard number $g^2$; if $B$ does not have CM, then
\[ \rho(A) = {{g+1} \choose 2} \leq 1+(g-1)^2\]
because $g \geq 4$. The case of a self-product of an elliptic curve being dealt with, we can assume $k \leq g/2$. Then, by Corollary \ref{lemma1} we have
\[ \rho(B^k) \leq \frac{1}{2}g(2k+1) \leq \frac{1}{2}g(g+1)\]
and the claim follows, since the equality 
\[ \frac{1}{2}g(g+1) \leq (g-1)^2 +1\]
holds for $g \geq 4$.\qed

\subsection{Proof of part (2)}
To start with observe that, if $r(A) \geq 3$, then
\[ \rho(A) \leq M_{r(A),g} \leq M_{3,g} < (g-2)^2 +4. \]
Therefore, we can assume $r(A) \leq 2$. Suppose that $A$ is an abelian variety with $r(A) =1$, i.e. $A \sim B^k$ with  $\dim B = b$ and $bk =g$. If $b=1$, then $B$ is an elliptic curve and we have two cases according to whether $B$ has complex multiplication. If $B$ does have complex multiplication, then $\rho(A) = g^2$ (the maximal Picard number), otherwise 
$\rho(A) = \frac{1}{2}g(g+1) < (g-2)^2+4$ (as $g \geq 7$).  If $b >1$, then $k \leq g/2$ and thus, by Corollary \ref{lemma1}, 
\[ \rho(B^k) \leq \frac{1}{2}g(2k+1) \leq \frac{1}{2}g(g+1) \leq (g-2)^2 +4\]
again because $g \geq 7$. The last remaining case is $r(A) = 2$, which we will divide into three steps.

\subsection*{Step 1}
We deal with abelian varieties of the form $E_1^n \times E_2^{g-n}$, where $E_1$ and $E_2$ are elliptic curves, and $1 \leq n \leq g-n$. If $n =1$, then, by Proposition \ref{picnum}
\[ \rho(E_1 \times E_2^{g-1}) = 1+ \rho(E_2^{g-1})\]
which equals $M_{2,g}$ if $E_2$ has complex multiplication, and $1 + \frac{1}{2}g(g-1)$ otherwise. In the CM case, we obtain the second largest attainable Picard number, in the non-CM case instead one sees that it is always the case that $1 + \frac{1}{2}g(g-1) \leq (g-2)^2 + 4$. Suppose now that $n \geq 2$: we have that $\rho(E_1^n \times E_2^{g-n}) \leq n^2 + (g-n)^2$, and we want to bound the right-hand side. The function
\[ f(x):= x^2 + (g-x)^2 \]
attains its maximum on the interval $\Omega=[2,g-2]$ at $x=2$ and $x=g-2$ with value $f(2)=f(g-2)=(g-2)^2+4$. This implies that $\rho(E_1^n \times E_2^{g-n}) \leq (g-2)^2+4$.

\subsection*{Step 2}
We now consider abelian varieties of the form $E^k \times A^l$, with $E$ an elliptic curve, $\dim A = a >1$, $k\geq 1$, $l\geq 1$ and $g = k+al$. Notice that, by Proposition \ref{picnum} and Lemma \ref{lemma1}, one has
\[ \rho(E^k \times A^l) \leq k^2 + \frac{1}{2}al(2l+1) = k^2 + \frac{1}{2}(g-k)(2l+1).\]
Consider the function
\[ f(x,y)= x^2 + \frac{1}{2}(g-x)(2y+1),\]
in the domain $\Omega := \lbrace (x,y) \in \R^2 \, \vert \, x \geq 1, \ y \geq 1, \ x+2y \leq g\rbrace$. We will prove that $f$ is bounded from above by $(g-2)^2 +4$ in $\Omega$.

By looking at the partials
\[ \frac{\partial f}{\partial x}(x,y) = 2x -y- \frac{1}{2}, \qquad \frac{\partial f}{\partial y}(x,y) = g-x,\]
we see that $f$ is increasing on the lines where $x$ is constant. Thus the maximum of $f$ in $\Omega$ will lie on the line $x + 2y = g$. Therefore, we have reduced ourselves to studying the function
\[ g(y) := f(g-2y,y) = (g-2y)^2 + 2y^2+ y\]
on $[1,(g-1)/2]$. Its maximum is at $y_{\text{max}}= 1$, with value 
\[ g(y_{\text{max}}) = (g-2)^2+3 < (g-2)^2 + 4.\]

\subsection*{Step 3}
The last case is that of products of the form $A^k \times B^l$, with $\dim A = a >1$, $\dim B = b >1$, $k  \geq l \geq 1$ and $g = ak+bl$. One has,
\begin{align*}
\rho(A^k \times B^l) &\leq \frac{1}{2}ak(2k+1) + \frac{1}{2}bl(2l+1) \\
& \leq \frac{1}{2} ak (2k+1) + \frac{1}{2}bl(2k+1) = \\
& = \frac{1}{2} g(2k+1) \leq \frac{1}{2}g(g-1) < (g-2)^2 + 4.
\end{align*}
\qed

\section{Structure of abelian varieties with large Picard number}\label{sec:Structure}
As an application of Theorem \ref{mainthm1}, we will now derive a structure result for abelian varieties of large Picard number (up to isogeny). Our starting point is the following result:
\begin{thm}[Exercise 5.6.10 of \cite{birkenhake-lange04}]\label{maxpic}
	Let $A$ be an abelian variety of dimension $g$. The following are equivalent 
	\begin{enumerate}
		\item $\rho(A) = g^2$;
		\item $A \sim E^g$, for some elliptic curve $E$ with complex multiplication;
		\item $A \cong E_1 \times \cdots \times E_g$, for some pairwise isogenous elliptic curves $E_1, \dots, E_g$ with complex multiplication.
	\end{enumerate}
\end{thm}

This result points out how the Picard number can force the structure of an algebraic variety to be in some sense rigid. Algebraic varieties with the maximum Picard number possible have shown to possess interesting arithmetic and geometric properties: for example, see \cite{shioda-mitani74} and \cite{shioda-inose77}, or \cite{beauville14} for a recent account. 

The aim of this section is to prove a similar statement for abelian varieties whose Picard number is the second or third largest attainable according to Theorem \ref{mainthm1}, namely $(g-1)^2+1$ or $(g-2)^2+4$. However, unlike in the case of maximal Picard number, one cannot expect a statement which is analogous to Theorem \ref{maxpic}(3). Already for $\rho(A)=(g-1)^2+1$, one can construct abelian varieties which are isogenous to $E_1^{g-1} \times E_2$, but which are not isomorphic to a product of elliptic curves.

The following result describes the structure of these abelian varieties up to isogeny. It should be noticed that, in contrast with Theorem \ref{maxpic}, the result depends on the dimension of the abelian varieties we consider: on one hand we need Theorem \ref{mainthm1}, and on the other we need to guarantee that the abelian varieties having such Picard numbers all belong to a unique isogeny class.

\begin{thm}\label{structurethm}
	Let $A$ be an abelian variety of dimension $g$.
	\begin{enumerate}
		\item Suppose $g\geq 5$. Then,
		\[ \rho(A) = (g-1)^2+1 \Longleftrightarrow A \sim E_1^{g-1} \times E_2, \]
		where $E_1$ has complex multiplication and $E_1$ and $E_2$ are not isogeneous.
			\item Suppose $g\geq 7$. Then,
			\[ \rho(A) = (g-2)^2+4 \Longleftrightarrow A \sim E_1^{g-2} \times E_2^2, \]
			where $E_1$ and $E_2$ both have complex multiplication but are not isogeneous. 
	\end{enumerate}
\end{thm}

\begin{proof}
	Recall that we have the following (strictly) increasing sequence of positive integers:
	\[ g=M_{g,g} < M_{g-1,g} < \cdots < M_{3,g} < M_{2,g} < M_{1,g}=g^2. \]
	Assume $\rho(A)=(g-1)^2+1 = M_{2,g}$ and $g \geq 5$. By definition of $M_{r,g}$ it follows that $r(A) \leq 2$; we claim that $r(A) =2$. Indeed, if $r(A)=1$, then necessarily $A \sim E^g$. But then either $\rho(A)=g^2$ if $E$ has CM (by Theorem \ref{maxpic}) or 
	$\rho(A)={{g+1} \choose 2} $ otherwise, either of which is a contradiction.
    Therefore $r(A)=2$ and (1) follows from Corollary \ref{lemmastr}.
	
	Now let $A$ have Picard number $\rho(A) = (g-2)^2+4$, and let $g \geq 7$. As $\rho(A) > M_{3,g}$, we deduce $r(A)\leq 2$. If $r(A)=1$ we again get a contradiction as above 
	(this time, we  also use $g \geq 7$), hence $r(A)=2$. In a similar fashion to the proof of the Theorem \ref{mainthm1}, we distinguish three cases:
	\begin{enumerate}
		\item[(a)] Let $A \sim A^k \times B^l$, with $\dim A >1$ and $\dim B >1$. Then, as we have seen in Step 3 of the proof of Theorem \ref{mainthm1}
		\[\rho(A) \leq \frac{1}{2} g(g-1) < (g-2)^2 +4,\] 
		which gives a contradiction.
		\item[(b)] Let $A \sim E^k \times B^l$, with $\dim B >1$ and $\dim E =1$. Then, as we have seen in Step 2 
		\[\rho(A) \leq (g-2)^2 +3 < (g-2)^2 +4,\]
		again a contradiction.
		\item[(c)] Let $A \sim E_1^n \times E_2^{g-n}$, for two elliptic curves $E_1$ and $E_2$. Then, the cases $n=1$ and $n=g-1$ can be discarded, by previous discussions.                  
		Therefore, let us suppose $2 \leq n \leq g-2$. We claim that $\rho(A) < (g-2)^2 + 4$, unless both $E_1$ and $E_2$ have complex multiplication. Indeed, if one of the factors
		 does not have complex multiplication then $\rho(A) \leq 3+(g-2)^2 < 4+(g-2)^2$. 			Therefore both $E_1$ and $E_2$ must have complex multiplication, and so $\rho(A)= \rho(E_1^n \times E_2^{g-n}) = n^2 + (g-n)^2$. 
			The maximum of this expression is achieved for $n=2$ or $n=g-2$, and this corresponds to a product $E_1^2 \times E_2^{g-2}$.
	\end{enumerate}
	We have thus shown that the only possible case is $A \sim E_1^2 \times E_2^{g-2}$, for two non-isogeneous elliptic curves $E_1$ and $E_2$ with complex multiplication, hence proving (2), and in this case the Picard number is  as stated. 
\end{proof}

Clearly, one can continue this analysis along arguments used above. However, one cannot,  in general,  expect to obtain a unique decomposition for a given Picard number. Already for  $\rho = (g-2)^2 +3$ there are two possible isogeny decompositions, namely:
\begin{enumerate}
	\item $E_1^{g-2} \times E_2^2$, $E_1$ being an elliptic curve with complex multiplication and $E_2$ not having complex multiplication;
	
	\item $E^{g-2} \times S$, $E$ being an elliptic curve with complex multiplication and $S$ being a simple abelian surface of type II (these do exist by results of Shimura \cite{shimura63}, see also the discussion in Section \ref{counting_picard_numbers}).
\end{enumerate}

\begin{rmk}
The Picard numbers $g^2$ and $(g-2)^2+4$ both lead to cases which have no complex moduli, whereas the intermediate Picard number $(g-1)^2+1$ leads 
to $1$-dimensional families. This is in striking contrast to the case of K3 surfaces where increasing the Picard number by one corresponds to a decrease in the number of moduli
by one. This is clear from the Torelli theorem for K3 surfaces. The difference lies in the fact that the Torelli theorem for K3 surfaces works with a weight $2$ Hodge structure, wheres abelian varieties  are governed by weight $1$ Hodge structures.
\end{rmk}

\section{Computing Picard numbers}\label{counting_picard_numbers}

In this section we approach the question how to compute the set $R_g$  of possible  Picard numbers of abelian varieties of a given dimension $g$. To this end, let us fix a positive integer $G$, such that we are interested in computing $R_G$. Because of the structure of $R_G$, we will in fact have to compute the sets $R_g$ for all $g \leq G$. In order to do this, we need to compute the Picard numbers of simple abelian varieties of dimension $g$, for every $g \leq G$. 

\subsection{Picard numbers of simple abelian varieties}
Let $g \geq 1$ be a fixed integer. If $X$ is a simple abelian variety of dimension $g$, its Picard number $\rho = \rho(X)$ must respect some divisibilty conditions \cite[Proposition 5.5.7]{birkenhake-lange04}, namely
\begin{enumerate}
	\item[($I$)] \textit{Type I}: $\rho \vert g$;
	\item[($II$)] \textit{Type II}: $\rho \in 3\IN$ and $\frac{2}{3}\rho \vert g$;
	\item[($III$)] \textit{Type III}: $2 \rho \vert g$;
	\item[($IV$)] \textit{Type IV}: $\rho \vert g$.
\end{enumerate}

For a fixed dimension $g$, we would like to understand which $\rho$ satisfying condition ($I$), ($II$), ($III$ or ($IV$) above can actually appear as the Picard number of a simple abelian variety of the corresponding type. For the notation used in the statement and in the proof of the next result, we refer the reader to Section \ref{sec:prelim} and to \cite[Ch.~9]{birkenhake-lange04} (in particular Section 9.6).

\begin{prop}\label{simple}
	Let $g$ be a fixed positive integer. For all positive integers $\rho$ that satisfy one of the conditions above, there exists a simple abelian variety $X$ of the corresponding type such that $\rho(X) =\rho$, unless we are in one of the five following exceptional cases:
	\begin{enumerate}
		\item $F$ is of type III, and $m := g/2e =1$;
		\item $F$ is of type III, $m := g/2e =2$, and there exists a totally positive element $\alpha \in K$ such that $N(T) = \alpha^2$ ($N$ being the reduced norm of $M_2(F)$ to $K$);
		\item $F$ is of type IV, $\sum_{\nu=1}^{e_0} r_\nu s_\nu =0$;
		\item $F$ is of type IV, $m:=g/d^2e_0=2$, $d=1$ and $r_\nu = s_\nu =1$ for all $\nu =1, \dots, e_0$;
		\item $F$ is of type IV, $m:=g/d^2e_0=1$, $d=2$ and $r_\nu = s_\nu =1$ for all $\nu =1, \dots, e_0$.
	\end{enumerate}
\end{prop}

\begin{proof}
	It is a theorem of Shimura that given an endomorphism structure $(F,\,',\iota)$ one has that a general member $(X,H,\iota)$ of the moduli space $\ka(\km,T)$ has the property $\End_\Q(X) = \iota(F)$, except in the cases above (for example see \cite{shimura63}, or \cite{birkenhake-lange04} for a modern approach). 

	In fact, under the assumption that our abelian variety $X$ be simple, one can show that these cases never occur:
	\begin{enumerate}
		\item $X$ is isogenous to a square $Y^2$, where $Y$ is an abelian variety of dimension $e_0$, contradicting the fact that $X$ is simple;
		\item same argument as above;
		\item $X$ is isogenous to $Y^{d^2m}$, where $Y$ is an abelian variety of dimension $e_0$, thus $d=m=1$, while $d=2$ because $F$ is a quaternion algebra over its center;
		\item $\End_\Q(X)$ contains a totally indefinite quaternion algebra $\tilde{F}$ over $K_0$ with $F=K \subset \tilde{F}$, so that $F = K \subset \tilde{F} \subset End_\Q(X) = F$, contradiction;
		\item as in (1) and (2).
	\end{enumerate}
	For details, consult \cite[Ch.~9, Ex.~9.10(1)--(5)]{birkenhake-lange04}, or see the original paper by Shimura \cite{shimura63}. By work of Gerritzen \cite{gerritzen71} and Albert \cite{albert34a,albert34,albert35}, this implies that given an involutorial division algebra $F$ of type I--IV outside of the five exceptional cases above, there exists a simple abelian variety whose endomorphism algebra is $F$. For a survey on these results, see \cite{oort88}. 
	
	We are now left to show that for any integer $\rho$ satisfying one of the conditions (I-IV) and outside of (1)-(5), we can actually construct an involutorial division algebra $F$ of the corresponding type, such that there exists an abelian variety $X$ with $\End_\Q(X) \cong F$ and $\rho(X) = \rho$. We will divide our analysis according to the type.
	
	\subsubsection{Type I}
	Let $\rho$ be a positive integer such that $\rho \vert g$. In this case, it is enough to construct a totally real number field $F$ of degree $\rho$ over $\Q$. However, given a finite abelian group $G$, it is always possible to construct a totally real number field $F$ such that $\Gal(F/\Q) \cong G$ as a subfield of a suitable cyclotomic field. This implies, in particular, that we can exhibit a totally real number field of degree $\rho$ over $\Q$.
	
	\subsubsection{Type II}
	It is enough to construct a totally indefinite quaternion algebra $F$ over a totally real number field $K$ of degree $e=[K:\Q]$ over $\Q$, such that $F$ is a division ring, as then $\rho = 3e$ satisfies the required condition. We are able to exhibit such an algebra for any $e \vert g$ by simply considering a quaternion algebra $F$ whose ramification is non-empty and disjoint from the archimedean place of $K$, i.e.
	\[ \emptyset = \text{Ram}(F) \cap \lbrace \text{archimedean places of $K$} \rbrace \quad \text{and} \quad \text{Ram}(F) \neq \emptyset. \]
	
	\subsubsection{Type III}
	In this situation, we aim at constructing a totally definite quaternion algebra $F$ over a totally real number field $K$ of degree $e=[K:\Q]$ over $\Q$, such that $F$ is in fact a division ring. It is enough to consider a quaternion algebra $F$ whose ramification locus is non-empty (this ensures the condition of being a division ring) and fulfills the condition 
	\[\emptyset \neq \text{Ram}(F) \supseteq \lbrace \text{archimedean places of $K$} \rbrace.\]
	
	\subsubsection{Type IV}
	We are left with the case corresponding to an involutorial division algebra $F$ of the second kind, whose center $K$ is a CM field with maximal real subfield $K_0$. In this case, it is enough to consider CM fields $K$ such that the degree $e_0=[K_0:\Q]$ of the maximal totally real subfield $K_0$ of $K$ ranges among all divisors of $g$ (i.e.~we are considering the case $d=1$). 
	
\end{proof}

\section{Additivity of the range of Picard numbers}\label{sec:additivity}

\subsection{Additivity}
As before we  denote by $R_g$ the set of realizable Picard numbers of $g$-dimensional abelian varieties, i.e.
\[ R_g := \big\lbrace \rho \, \vert \, \text{$\exists X$ abelian variety, $\dim X= g $, $\rho(X) =\rho$} \big\rbrace. \]
Conventionally, let us set $R_0 = \lbrace 0 \rbrace$. The main result of this section is 
\begin{prop}\label{prop:additivity} 
For any integers $g,h \geq 0$ we have an inclusion 
\[ R_g + R_h := \lbrace x+y \ \vert \ x \in R_g, \ y\in R_h \rbrace \subset R_{g+h}. \]
\end{prop}

Clearly, the idea is to use the additivity of Picard numbers for products of non-isogeneous abelian varieties, as proven in Corollary \ref{picnum}.
In order to be able to us this we need that for any $k\in R_g$ we find countably many abelian varieties of dimension $g$ and Picard number $k$ in different isogeny classes.
The case of elliptic curves illustrates how this can be proved:    
Suppose $E$ and $E'$ are two elliptic curves, and let $F:=\End_\Q (E)$ and $\End_\Q (E')$ be their endormorphism algebras. If $E$ and $E'$ have CM, then they are mutually isogenous if and only if $F \cong F'$. However, if they don't have CM, then $F \cong F' \cong \Q$ but $E$ and $E'$ are not necessarily isogenous. However, by removing CM elliptic curves from $\km_{1,1}$, we get an uncountable set of isomorphism classes of elliptic curves. Since each isogeny class consists of countably many (isomorphism classes of) elliptic curves, we must have infinitely many isogeny classes of non-CM elliptic curves. We first note
\begin{prop}
	Let $X$ and $X'$ be two simple abelian varieties of dimension $g$, and let $F$ and $F'$ be the corresponding endomorphisms algebras. Then, if $X$ is isogenous to $X'$, then $F \cong F'$.
\end{prop}

\begin{proof}
	The isomorphism $\iota: \End_\Q{X'} \longra \End_\Q{X}$ is defined by sending $\alpha  \longmapsto \psi \circ \alpha \circ \phi$, where $\phi: X \longra X'$ is an isogeny, and $\psi: X' \longra X$ is the unique isogeny such that $\psi \circ \phi= e_X$ and $\phi \circ \psi= e_Y$, $e$ being the exponent of $\phi$ (or  $\psi$ respectively). Surjectivity and injectivity of $\iota$ follow from the fact that multiplication maps are invertible in the endomorphism $\Q$-algebra.
\end{proof}
The key result of this section is the 
\begin{prop}\label{prop:infinitnelymany}
Given an integer $k \in R_g$, then there exist at least countably many isogeny classes of abelian varietiss of dimension $g$ and Picard number $k$.  
\end{prop}
\begin{proof}
Assume that $k \in R_g$ and that this integer is realized by an abelian variety $A$ of  dimension $g$. We first assume that $F:=\End_\Q(X) \not\cong \Q$.
We will now go through the various types of endomorphism algebras and start with type I and $[F:\Q] >1$. 

For type I, the endomorphism algebra is a totally real number field $F$, and the Picard number of a simple abelian variety with such an endomorphism algebra is the degree $e:=[F:\Q]$ of $F$. One can find infinitely many totally real number fields of a fixed degree $e>1$. By taking the corresponding (infinitely many) abelian varieties, we are done. 

For type II-III, the endomorphism algebra is a quaternion algebra $F$ over a totally real number field $K$. In this case, the Picard number is $\rho=3e$ for type II and $\rho = e$ for type III, where $e:=[K:\Q]$. Since by the previous argument there exist infinitely many totally real number fields of degree $e$, it follows that we can find infinitely many (totally definite or totally indefinite) quaternion algebras. By considering the corresponding abelian varieties, we have shown the claim for types II-III.

For type IV, the endomorphism algebra $F$ has degree $[F:K]=d^2$ over a CM field $K$, $e:=[K:\Q]$. If $K_0$ is the totally real subfield of $K$, of degree $e_0:=[K_0:\Q]$, the Picard number is $\rho= e_0d^2$. Similarly to the argument in the proof of Proposition \ref{simple}, we can restrict ourselves to consider abelian varieties whose endomorphism algebra satisfies the condition $d=1$: under this assumption, $\rho=e_0$. Since there are infinitely many totally real number fields, we find infinitely many CM fields (by just adding the imaginary unit), and the same argument as in the previous cases shows the statement for type IV abelian varieties.

This leaves us with the general situation where $F=\End_\Q(X) \cong \Q$. The ppav of this type are given by removing from $\mathcal A_g$ a countable union of proper Shimura varieties.
Since $\mathcal A_g$ has positive dimension and since the set of ppav isogeneous to a given abelian variety is countable, the claim follows.
\end{proof}

\begin{proof}[Proof of Proposition \ref{prop:additivity}]
This now follows immediately from Proposition \ref{prop:infinitnelymany} and the additivity proved in Corollary \ref{picnum}.
\end{proof}

\subsection{Computing $R_g$}\label{algorithm}

Our final aim is to find all realizable Picard number of abelian varieties of a given dimension. We can use Proposition \ref{prop:additivity} to easily show that some of the lower one indeed occur.

\begin{prop}\label{products}
	Given $g \geq 2$, consider the set $R_g$ of Picard numbers of abelian varieties of dimension $g$. Then, $\lbrace 1, \dots, 2g \rbrace \subset R_g$. 
\end{prop}

\begin{proof}
	As $R_2 =  \lbrace 1, \dots, 4 \rbrace$ is complete and $R_3 \supset \lbrace 1, \dots, 6 \rbrace$, the result easily follows by induction on $g$.
\end{proof}

\begin{rmk}
	In fact, it is not hard to prove that all Picard numbers $\rho$ satisfying the inequality $g \leq \rho \leq 2g$ are attained by products of elliptic curves.
\end{rmk}

This now allows us to formulate an algorithm  which computes the ranges $R_g$ inductively.

\begin{itemize}
	\item set $R_1 = \lbrace 1 \rbrace$ and $R_2 = \lbrace 1,2,3,4 \rbrace$;
	\item for all $g$ in the range $3 \leq g \leq G$, we compute $R_g$ as follows:
	\begin{enumerate}
		\item[($\rm{i}$)] by Proposition \ref{products}, $R_g \supset \lbrace 1, \dots, 2g-1\}$ (in particular, all Picard numbers of simple abelian varieties of dimension $g$ are in this range);
		\item[($\rm{ii}$)] compute all possible Picard numbers of self-product abelian varieties $A^k$, where $\dim A = g/k$;
		\item[($\rm{iiii}$)] for every pair $(g_1,g_2)$ of positive integers such that $g_1 + g_2 =g$, compute $R_{g_1} + R_{g_2}$;
		\item[($\rm{iv}$)] assemble everything in light of 
		\[ 
		R_g = \bigcup_{k \vert g} \big\lbrace \rho(A^k) \, \vert \, \text{$A$ simple, $\dim A=g/k$} \rbrace \cup \bigcup_{1 \leq n \leq g-1} \big(R_{n}+R_{g-n}\big).
		\]
	\end{enumerate}  	
\end{itemize}

\section{Asymptotic behaviour of Picard numbers of abelian varieties}\label{sec:asymptotics}

\subsection{Asymptotic completeness of Picard numbers}
In the course of this note, we have shown that for every $g \geq 3$ the set $R_g$ of Picard numbers of $g$-dimensional abelian varieties is not complete, or in other words that $\# R_g < g^2$. The ratio $\delta_g := \#R_g/g^2$ is the \emph{density} of $R_g$ in $[1,g^2]\cap\IN$, and it describes how many admissible Picard numbers (according to Lefschetz Theorem of $(1,1)$-classes) can actually be attained. One may wonder about the \textit{asymptotic density} of Picard numbers of abelian varieties: this is the quantity defined as
\[ \delta := \lim_{g \ra +\infty} \delta_g. \]
We now show that the Picard numbers of abelian varieties are asymptotically complete, namely that $\delta =1$, contrary to the fact that $\delta_g <1$ for every $g\geq 3$.

\begin{thm}[Asymptotic completeness]\label{density}
	The sets of Picard numbers of abelian varieties are asymptotically dense, i.e.~$\delta=1$.
\end{thm}

The proof relies on Lagrange's four-square theorem and the following lemma, whose proof follows readily from the additivity of the Picard number.

\begin{lemma}\label{lombardo}
	Suppose $g \geq 1$ and $1 \leq n \leq g^2$, where $g$ and $n$ are two integers. Assume that there exist positive integers $n_1, \dots, n_k$ such that 
	\[ n-1 = n_1^2 + \cdots + n_k^2 \qquad \text{and} \qquad n_1 + \cdots + n_k \leq g-1. \]
	Then, there exists a $g$-dimensional abelian variety $X$ with $\rho(X) =n$.
\end{lemma}

\begin{proof}[Proof of Proposition \ref{density}]
	Let $n_1$ be the largest positive integer such that
	\[ n_1^2 \leq n-1 < (n_1+1)^2. \]
	Then,
	\[ 0 \leq n-1-n_1^2 < (n_1+1)^2 - n_1^2 = 2n_1 +1, \]
	from which it follows that 
	\[ 0 \leq n-1-n_1^2 \leq 2n_1 \leq 2 \sqrt{n-1} < 2 \sqrt{n} \leq 2g. \]
	Lagrange's four-square theorem implies that
	\[ m:=n-1-n_1^2 = n_2^2 + n_3^2 + n_4^2 + n_5^2, \]
	for some $n_2, n_3, n_4, n_5 \in \IN$. We will now show that $n_1 + \cdots + n_5 < g$ for $g \gg 0$: indeed, by looking at the power means of $n_2, \dots, n_5$ one has that
	\[ \frac{n_2 + \cdots + n_5}{4} \leq \sqrt{\frac{n_2^2 + n_3^2 + n_4^2 + n_5^2}{4}} = \frac{1}{2} \sqrt{m} \leq \frac{1}{2} \sqrt{2g}. \]
	Therefore, 
	\[ n_1 + \cdots + n_5 \leq \sqrt{n-1} + 2 \sqrt{2g},\]
	and the right-hand side is strictly smaller than $g$ if and only if
	\[ n < g^2 + 8g+1 - 4\sqrt{2}g^{3/2} =: b_g.\]
	This implies that all Picard numbers in the range $[1, b_g)$ indeed occur, 
	by virtue of the lemma above. Hence, we have that asymptotically $\# R_g \geq b_g -1 = g^2 + 8g - 4\sqrt{2}g^{3/2}$, and thus $\delta =1$.
\end{proof}

\subsection{Distribution of large Picard numbers}
We are interested in describing the distribution of large Picard numbers within $[1,g^2] \cap \IN$. As we have already observed, for every $g \geq 1$, the set $R_g$ has the following structure
\[ 
R_g = \bigcup_{k \vert g} \big\lbrace \rho(A^k) \, \vert \, \text{$A$ simple, $\dim A=g/k$} \rbrace \cup \bigcup_{1 \leq n \leq g-1} \big(R_{n}+R_{g-n}\big).
\]

The proof of Theorem \ref{mainthm1}(1) shows that all Picard numbers of abelian varieties of dimension $g$ that are isogenous to a self-product of a simple abelian variety are bounded by $\frac{1}{2}g(g+1)$, unless we are considering the $g$-fold product of a CM elliptic curve, in which case the maximal Picard number is attained. 

In order to begin our analysis, we need to specify what we mean by "large Picard numbers". First of all, we will require large Picard numbers to satisfy the inequality $\rho > g(g+1)/2$. In particular, this implies that we need not concern ourselves with those abelian varieties whose isogeny decomposition only has one factor. Therefore, we can focus on the following subset of $R_g$:
\[ \bigcup_{1 \leq n \leq g-1} \big(R_{n}+R_{g-n}\big). \]

Now we need a little bit of notation. Let us set
\[ R_{g,n} := \big\lbrace (g-n)^2 + x \, \vert \, x \in R_n \big\rbrace. \]
In other words, $R_{g,n}$ is the subset of $R_g$ obtained by translating $R_n$ to the right by $(g-n)^2$ inside $\IN$. Notice that given $g,k,n \in \IN$, one has by Proposition \ref{prop:additivity} that 
\[ R_{g,k} + R_n \subset R_{g+n,k+n}.  \] 

As we want to consider large Picard numbers only, we will be concerned only with some of the $R_{g,s}$'s, namely those for which the inequality
\begin{equation}\label{equ:notsimple1}
 \frac{1}{2}g(g+1) \leq (g-s)^2 +1 
 \end{equation}
holds ($1 \leq s \leq g$), which implies that the abelian varieties we are considering are not self-products of simple abelian varieties. This in particular implies that $g \geq 4$ (because $s \geq 1$) and
\begin{equation*}
s \leq \frac{2g-\sqrt{2(g^2+g-2)}}{2}.
\end{equation*}

Finally, let us look at the mutual interaction of the $R_{g,s}$'s. For a fixed $g$, there might exist positive integers $a$ and $b$ such that $R_{g,a} \cap R_{g,b} \neq \emptyset$. However, if $a$ and $b$ are distinct and small enough with respect to $g$, then $R_{g,a} \cap R_{g,b} = \emptyset$. Indeed, for a fixed $g$, the inequality 
\begin{equation}
	[g-(s+1)]^2+(s+1)^2 < (g-s)^2 +1
\end{equation}
holds for all positive integers $s < -2 + \sqrt{2g+3}$. Hence $R_{g,a} \cap R_{g,b} = \emptyset$ for $a$ and $b$ in this range. 

We are now able to define precisely what "large Picard number" means. We will say that a Picard number $\rho$ is \emph{large} if $\rho \in R_{g,s}$ with $s$ satisfying conditions ($1$) and ($2$). This condition can be made explicit:

\begin{prop}\label{lemmaasymptotic}
	Let $g \geq 4$. Then, $\rho \in R_{g,s}$ is a large Picard number if and only if
	\[ s \leq \min \Big\{ \frac{2g-\sqrt{2(g^2+g-2)}}{2}, -1 + \sqrt{2g+3} \Big\}.\]
\end{prop}

The theorem we are going to discuss next describes the distribution of the large Picard numbers inside $[1,g^2] \cap \IN$. The argument we use is inductive and its initial step is a proof of an asymptotic version of Theorem \ref{mainthm1}(2) starting from an asymptotic version of Theorem \ref{mainthm1}(1).

\begin{thm}[Distribution of large Picard numbers]\label{asymptoticstr}
        For every positive integer $\ell$ there exists a genus $g_\ell$ such that for all $g\geq g_\ell$ large Picard numbers in $R_g$ are distributed as follows:  
	\[  \boxed{R_{g,\ell}} \qquad  \cdots  \qquad \boxed{R_{g,4}} \qquad \boxed{R_{g,3}} \qquad \boxed{R_{g,2}} \qquad \bullet^{(g-1)^2+1} \qquad \bullet^{g^2}.  \]
	In other words, for all $g \geq g_\ell$, we have that
		\[ [(g-\ell)^2 +1, g^2] \cap R_g = R_{g,\ell} \sqcup R_{g,\ell-1} \sqcup \dots \sqcup R_{g,2} \sqcup R_{g,1} \sqcup R_{g,0}. 	\]
\end{thm}
\begin{rmk}
In particular this shows that, as $g \to \infty$ more and more gaps arise in $R_g$ as we go down from the maximum Picard number $g^2$. 
\end{rmk}

\begin{proof}
	We will give a proof by induction. To start with, note that we can (and will) always assume that $g$ is large enough, so that there is no overlapping between the sets $R_{g,n}$ that we wish to consider. 
	
Let us start with a pair $(\ell,g_\ell)$ such that the following holds: there is no abelian variety of dimension $g\geq g_\ell$ and Picard number 
$\rho$ such that $(g-t)^2 +t^2 < \rho < (g-t+1)^2 +1$ for $2 \leq t \leq \ell$ or $(g-1)^2+1 < \rho < g^2$. Then, the claim is that we can find  $g_{\ell+1}\geq g_\ell$ such that there is no abelian variety of dimension $g \geq g_{\ell+1}$ and Picard number $\rho$ such that $[g-(\ell+1)]^2+(\ell+1)^2 < \rho < (g-\ell)^2+1$.

	As the start of the induction, we will now show how to recover the second part of Theorem \ref{mainthm1} from the first one, at least asymptotically (we will not be able to get any bound on $g$, but of course it is always possible to do so). Suppose that for all $g \geq g_1$ there is no abelian variety $X$ of dimension $g$ and Picard number $(g-1)^2 + 1 < \rho(X) < g^2$ (notice that in light of Theorem \ref{mainthm1}, we can choose $g_1=4$). We will now prove that there exists $g_2$, $g_2 \geq g_1$, such that for every $g \geq g_2$ there is no abelian variety $Y$ of dimension $g$ and Picard number $(g-2)^2 + 4 < \rho(Y) < (g-1)^2 +1$.
	
	Let $g_2$ be such that for all $g \geq g_2$ conditions ($1$) and ($2$) above are satisfied (i.e. the Picard numbers $\rho$ in the range $(g-2)^2 + 4 < \rho< (g-1)^2 +1$ are large according to our definition). Suppose $Y$ is an abelian variety whose Picard number contradicts the statement we want to prove, namely $(g-2)^2 + 4 < \rho(Y) < (g-1)^2 +1$. Then, as $\rho(Y)$ is large, $Y$ is isogenous to a product of two abelian varieties, i.e.~$Y \sim A_n \times A_{g-n}$, where $n \leq g-n$ and $\Hom(A_n, A_{g-n})=0$, (here the subscripts indicate the dimension). Since $\rho > (g-2)^2 + 4$, we have that $n=1$ necessarily. Therefore $Y \sim E \times A_{g-1}$, where $E$ is an elliptic curve and $\Hom(E, A_{g-1})=0$. As $\rho(Y) = 1 + \rho(A_{g-1})$, we readily see that
	\[ (g-2)^2 + 1 < \rho(A_{g-1}) < (g-1)^2,  \]
	a contradiction. This is first step of the induction.
	
Now, let us assume that there exists $g_\ell$ such that for all $g \geq g_\ell$ there is no abelian variety $X$ of dimension $g$ and Picard number in the following ranges:
\begin{table}[h]
		\centering
		\begin{tabular}{lllll}
			&(1)  &$(g-1)^2 +1 < \rho(X) < g^2$; & &\\
			&(2)  &$(g-2)^2 + 4 < \rho(X) < (g-1)^2 +1$;  &  &  \\
			&  & $\qquad \qquad \qquad \vdots$ &  & \\
			&($\ell$) &$(g-\ell)^2 + \ell^2 < \rho(X) < [g-(\ell-1)]^2 + 1$.  & & 
		\end{tabular}
	\end{table}
	
We claim that there exists $g_{\ell+1}$ such that for all $g \geq g_{\ell+1}$ there is no abelian variety $Y$ of dimension $g$ and Picard number satisfying
\begin{table}[h]
		\centering
		\begin{tabular}{lllll}
			&($\ell+1$) &$[g-(\ell+1)]^2 + (\ell+1)^2 < \rho(Y) < [g-\ell]^2 + 1$.  & & 
		\end{tabular}
	\end{table}

Again, let us let $g$ grow bigger so that the Picard numbers we wish to consider can only be realized by abelian varieties that are not a self-product of a simple abelian variety. By contradition, let $Y$ be an abelian variety that contradicts the statement we want to prove. Then $Y \sim A_n \times A_{g-n}$, where $n \leq g-n$ and $\Hom(A_n, A_{g-n})=0$. It is straightforward to see that $n \leq \ell$, as $\rho(Y) >[g-(\ell+1)]^2 + (\ell+1)^2$. By additivity of the Picard number
\[ \rho(Y) = \rho(A_n \times A_{g-n}) = \underbrace{\rho(A_n)}_{\rho_n}+\underbrace{\rho(A_{g-n})}_{\rho_{g-n}}.\]
As $\rho(Y) >[g-(\ell+1)]^2 + (\ell+1)^2$, we see that
\begin{align*}
		\rho_{g-n} &> [g-(\ell+1)]^2 + (\ell+1)^2 - \rho_n > [g-(\ell+1)]^2 + (\ell+1)^2 -n^2\\
				  &> [(g-n)-(\ell-n+1)]^2 +(\ell+1)^2 - n^2 \\
				  &> [(g-n)-(\ell-n+1)]^2 +(\ell-n+1)^2.
\end{align*}
Similarly, 
\[ \rho_{g-n} < (g-\ell)^2 +1 - \rho_n \leq (g-\ell)^2 = [(g-n)-(\ell-n)]^2 < [(g-n)-(\ell-n)]^2 +1,\]
and summing up we have shown that 
\[ [(g-n)-(\ell-n+1)]^2 +(\ell-n+1)^2 < \rho(A_{g-n}) <[(g-n)-(\ell-n)]^2 +1,\]
which contradicts the $(\ell-n+1)$-st condition above.
	
\end{proof}

As a striking consequence of Theorem \ref{asymptoticstr}, we get the following structure theorem for abelian varieties of large Picard number up to isogeny, which generalizes the results in Section \ref{sec:Structure}. As we had already noticed in Section \ref{sec:Structure}, we cannot expect a structure theorem up to isomorphism, hence this is the strongest result we could hope for.
\begin{cor}[Structure theorem for abelian varieties of large Picard number]\label{cor:morestructure}
	For every positive integer $\ell$ there exists a genus $g_\ell$ such that for all $g\geq g_\ell$ the following are equivalent:
	\begin{enumerate}
		\item $\rho(X) \in R_{g,n}$ for some $n \leq \ell$;
		\item $X \sim E_{g-n} \times A_n$, where $E$ is an elliptic curve with complex multiplication, $A_n$ is an abelian variety of dimension $n$, and $\Hom(E,A_n)=0$.
	\end{enumerate}
\end{cor}

\begin{proof}
		Let us set $g_\ell$ as in the proof of Theorem \ref{asymptoticstr}, and let $X$ be an abelian variety of Picard number $\rho(X) \in R_{g,n}$ for some $n \leq \ell$. By means of the Poincar\'e reducibility theorem, we can write $X \sim E^{t} \times A_{g-t}$, where $E$ is an elliptic curve with complex multiplication, $A_{g-t}$ is an abelian variety of dimension $g-t$, $\Hom(E,A_{g-t})=0$, and $t$ is the largest integer appearing as exponent of an elliptic curve with complex multiplication in the isogeny decomposition of $X$. Let us now set for simplicity $t=g-m$, so that $X \sim E^{g-m} \times A_m$. In particular, it follows that $\rho(X) \in R_{g,m}$. However, by Theorem \ref{structurethm}, $R_{g,n}$ cannot intersect $R_{g,m}$ unless $n=m$, from which the statement follows.
\end{proof}

\section{Abelian varieties defined over number fields} \label{sec:numbefields}

In this section we show that every realizable Picard number $\rho \in R_g$ can be obtained by an abelian variety defined over a number field.

\begin{thm}\label{teo:endonumberfield}
Let  $(X,\lambda)$ be a polarized complex abelian variety, let $D = \End^0(X)$ be the endomorphism algebra of $X$ and let $*$ be the Rosati involution on~$D$. Then there exists a polarized abelian variety over 
$\overline{\mathbb Q}$, or equivalently over a number field,  which has the same endomorphism algebra with involution $(D,*)$.
\end{thm} 
\begin{proof}[Proof (Ben Moonen)]
To prove the assertion, choose a $\mathbb Q$-subalgebra $R \subset \mathbb C$ of finite type and a polarized abelian scheme $(Y,\mu)$ 
over $S:=\Spec(R)$ such that $(Y,\mu) \otimes_R \mathbb C$ is isomorphic to $(X,\lambda)$ and such that all endomorphisms of~$X$ are defined over~$R$, in the sense that the natural map 
\[ \End^0(Y/R) \longra \End^0(X) \]
is an isomorphism. The existence of such a model follows from \cite[Proposition~(8.9.1)]{EGA_IV_4} together with the fact that $\End(X)$ is a finitely generated algebra (in fact, it is even finitely generated as an abelian group).  By construction, if $\eta$ is the generic point of~$S$, we have $\End^0(Y_\eta) \cong D$ as algebras with involution. If $s$ is a point of~$S$, we have a specialization homomorphism $i_s\colon \End^0(Y_\eta) \hookrightarrow \End^0(Y_s)$, and we are done if we can find a closed point~$s$ for which $i_s$ is an isomorphism.

Let $\ell$ be a prime number. For $s$ a point of~$S$, let $T_{\ell}(s):=T_\ell(Y_s)$ denote the $\ell$-adic Tate module of~$Y_s$, and let
\[\rho_s \colon \Gal\bigl(\overline{\kappa(s)}/\kappa(s)\bigr) \to {\rm GL}(V_{\ell}(s))\]
denote the Galois representation on $V_{\ell}(s) = T_{\ell}(s) \otimes_\Z \mathbb Q_\ell$. By a result of Faltings \cite[Theorem 1]{faltings84}, $\End^0(Y_s) \otimes \mathbb Q_\ell$ is the endomorphism algebra of~$V_{\ell}(s)$ as a Galois representation, namely
\[ \End^0(Y_s) \otimes \mathbb Q_\ell \ {\cong} \ \End_{\Gal(\overline{\kappa(s)}/\kappa(s))}\big(V_\ell(s) \big). \]
For $s \in S$, the image of~$\rho_s$ may be identified with a subgroup of $\im(\rho_\eta)$; the subgroup we obtain is independent of choices only up to conjugacy. By a result of Serre \cite{serre_to_ribet} (see also \cite[Proposition 1.3]{noot95}), there exist closed points $s \in S$ for which $\im(\rho_s) = \im(\rho_\eta)$, and for all such points the specialization map~$i_s$ on endomorphism algebras is an isomorphism (see \cite[Corollary 1.5]{noot95}). 
\end{proof}

As the Picard number only depends on $(D,*)$ this immediately  implies
\begin{cor}
Every realizabe Picard number $\rho \in R_g$ can be obtained by an abelian variety defined over a number field.
\end{cor}

\bibliographystyle{plain}
\bibliography{bib}{}

\end{document}